\documentclass[12pt]{amsart}                                    
\usepackage{amsmath,amssymb,amsthm}
\usepackage{comment}
\usepackage[usenames]{color}
\usepackage[all,cmtip]{xy}
\usepackage{tikz}\usetikzlibrary{cd}
\usepackage{hyperref}
\usepackage{url}
\usepackage{float}\restylefloat{figure}

\newtheorem{thm}{Theorem}

\newtheorem{proposition}[thm]{Proposition}

\newtheorem{conj}[thm]{Conjecture}

\newtheorem{df}[thm]{Definition}

\newtheorem{notat}[thm]{Notation}
\newtheorem{proviso}[thm]{Proviso}
\newtheorem{rmk}[thm]{Remark}
\newtheorem{term}[thm]{Terminology}

\newenvironment{ls}{\begin{itemize}}{\end{itemize}}

\newcommand{\scr}[1]{\ensuremath{\mathcal {#1}}}

\renewcommand{\phi}{\varphi}
\newcommand{\eps}{\varepsilon}

\newcommand{\notarrow}{\kern .42em\not\kern -.42em\longrightarrow}

\newcommand{\x}{\ensuremath{\otimes}}
\newcommand{\w}[3]{\ensuremath{#1:#2\to#3}}

\renewcommand{\circ}{;}

\newcommand{\ax}{\ensuremath{\alpha^\otimes}}
\newcommand{\ap}{\ensuremath{\alpha^{\oplus}}}
\newcommand{\C}{\ensuremath{\scr C}}
\newcommand{\gp}{\ensuremath{\gamma^{\oplus}}}
\newcommand{\gx}{\ensuremath{\gamma^\otimes}}
\newcommand{\Hom}{\ensuremath{\text{Hom\,}}}
\newcommand{\IS}[2]{\ensuremath%
  {\text{Intermediate step}~#1\text{ in simplification of~#2}}}
\newcommand{\op}{\,\oplus\,}

\begin{document}

\title{Braided Distributivity}
\thanks{Part of this work was done while the first author was a
visiting researcher and the second author a principal researcher in the Quantum Architectures and Computing (QuArC) group at Microsoft Research.}
\author{Andreas Blass}
\address{Mathematics Department\\
University of Michigan\\
Ann Arbor, MI 48109--1043, U.S.A.}
\email{ablass@umich.edu}
\author{Yuri Gurevich}
\address{Computer Science and Engineering\\
University of Michigan\\
Ann Arbor, MI  48109-2121, U.S.A}
\email{gurevich@umich.edu}

\begin{abstract}
\maketitle

In category-theoretic models for the anyon systems proposed for
topological quantum computing, the essential ingredients are two
monoidal structures, $\oplus$ and $\otimes$. The former is
symmetric but the latter is only braided, and $\otimes$ is
required to distribute over $\oplus$. What are the appropriate
coherence conditions for the distributivity isomorphisms? We
came to this question working on a simplification of the
category-theoretical foundation of topological quantum
computing, which is the intended application of the research
reported here.

This question was answered by Laplaza when both monoidal
structures are symmetric, but topological quantum computation
depends crucially on $\otimes$ being only braided, not
symmetric. We propose coherence conditions for distributivity in
this situation, and we prove that our conditions are\\
(a) strong enough to imply Laplaza's when the latter are
suitably formulated, and\\
(b) weak enough to hold when --- as in
the categories used to model anyons --- the additive structure
is that of an abelian category and the braided $\otimes$ is
additive.\\ Working on these results, we found a new
redundancy in Laplaza's conditions.
\end{abstract}

\section{Introduction} 
\label{sec:intro}

Although this paper is about pure category theory, its origin,
motivation, and intended use lie in the application of
categories to topological quantum computing, specifically to the
description of nonabelian anyons and the calculation of their
properties. Traditionally that application employs modular
tensor categories, which involve a great deal of
category-theoretic structure (abelian categories with braided
monoidal structure, duality, ribbon structure, and more); the
definition is given in detail in \cite{PP} and with some
emendations in \cite{G225}.  The book \cite{BK} provides an
extensive treatment of the theory of these categories.

For the purposes of quantum computation, the most important
information about these categories is the braiding structure.
When computing this braiding structure in some specific cases,
like the case of Fibonacci anyons \cite[\S8.5 of the published
version, \S5 on the arXiv]{G225},  we found that only a small
part of the modular tensor category structure is really used in
these computations.  In particular, the computations can be done
in a framework where the only morphisms are isomorphisms; that
is, the categories can all be taken to be groupoids.  Of course,
they will no longer be abelian categories; an abelian category
is a groupoid if and only if it is equivalent to the category
with one object and one morphism. The groupoid framework is
arguably simpler. In addition, it allows a presentation that
looks less like category theory and more like universal algebra,
and thus may be accessible for a broader audience. We have
prepared a presentation \cite{G239} of this framework, and the
present paper plays an important technical role there.

In modular tensor categories, there are two monoidal structures,
one (written $\oplus$) coming from the assumption that the
category is abelian, and one (written $\otimes$) that is assumed
separately. (Both are subject to some additional assumptions
that need not concern us here.)  In our new framework, we do not
have an abelian category, so both $\oplus$ and $\otimes$ need to
be assumed individually.  The former is a symmetric monoidal
structure; the latter is only a braided monoidal structure.
Furthermore, $\otimes$ must distribute over $\oplus$ (a
requirement that, in abelian categories, would follow from
assuming that $\otimes$ is an additive bifunctor), and all the
relevant isomorphisms (associativity, commutativity, unit, and
distributivity) must satisfy suitable requirements.

What, exactly, are the suitable requirements? For a symmetric
monoidal structure, like $\oplus$, the appropriate requirements
were found by Mac Lane \cite{MacLane1963}and subsequently
simplified by Kelly \cite{Kelly}. For a braided monoidal
structure, the requirements were supplied by Joyal and Street
when they introduced the notion of braided structure in
\cite{JS1,JS2}.

Finally, for distributivity, Laplaza \cite{Laplaza} has found
the requirements that should be satisfied by distributivity
isomorphisms (in fact, even by distributivity monomorphisms) for
a pair of symmetric monoidal structures. Of course, any
universal statement about monomorphisms implies the same
statement for isomorphisms. We are interested in the
isomorphism-bound case of Laplaza's study where distributivity
is assumed to be given by isomorphisms, not merely by
monomorphisms.

\begin{proviso}\label{prov:laplaza}\rm
Below, by default, speaking about Laplaza's coherence conditions, we mean the isomorphism-bound versions of his conditions.
\end{proviso}

Laplaza's coherence conditions do not give us exactly what we
need because he assumes that both $\oplus$ and $\otimes$ are
symmetric, whereas in our situation $\oplus$ is symmetric but
$\otimes$ is only braided.  The difference is important because
Laplaza's conditions (specifically, his condition II) require
the left distributivity $A\otimes(B\oplus C)\cong (A\otimes
B)\oplus(A\otimes C)$ and the right distributivity $(B\oplus
C)\otimes A\cong(B\otimes A)\oplus(C\otimes A)$ to be related
via the relevant commutativity isomorphisms. In the braided
case, there are two commutativity isomorphisms $X\otimes Y\cong
Y\otimes X$, so we need to either choose one to prefer or make
sure they transform left distributivity to right distributivity
in the same way.

Also, Laplaza lists 24 coherence conditions, but then proves
that many of them are redundant, thus providing a considerably
reduced but still sufficient list of coherence conditions.


We found one new redundancy in Laplaza's coherence conditions (in his original context and in the isomorphism-bound version of his context),
namely that condition~XVIII implies condition~XVII; see
Proposition~\ref{redundancy} in \S\ref{sec:axioms} below. This
allows us to slightly improve Laplaza's reduced list of
coherence conditions.

We need to check whether all these redundancies, including the
new one, apply in our context or whether some of them depend on
the symmetry (rather than mere braiding) of $\otimes$. It turns
out that they do apply.

The goal of this paper is to present carefully the appropriate
requirements for the situation of a braided monoidal structure
$\otimes$ distributing over a symmetric monoidal structure
$\oplus$.

\begin{term}\rm
In all the preceding cases, the requirements are called
coherence conditions and come with an associated coherence
theorem. But we do not provide an associated coherence theorem;
accordingly we call our requirements axioms. Throughout the
paper, we use the term \emph{requirement} as more general than
\emph{coherence condition} or \emph{axiom}, so that coherence
conditions are requirements and so are axioms.
\end{term}

In \S\ref{sec:axioms} we present our notion of a category equipped with a symmetric monoidal structure $\oplus$ and a braided monoidal structure $\otimes$ with $\otimes$ distributing over $\oplus$. For brevity, we call such categories BD categories; here ``BD" alludes to braided distributivity. The appropriate requirements for the distributivity isomorphisms are called BD axioms.

In \S\ref{sec:cons} we show that the BD axioms are strong enough by deducing from them all of Laplaza's coherence conditions appropriately formulated for the braided situation.  Additional evidence of their strength is given in \cite{G239}, where we show how they support the computation of associativity and braiding isomorphisms for particular anyons.

Finally, we show in \S\ref{sec:abel} that the BD axioms are not
too strong by showing that they hold when --- as in the
categories used to model anyons --- the additive structure is
that of an abelian category and the multiplication bifunctor of
the braided monoidal structure is additive.

Even apart from the motivation from topological quantum
computing, the determination of an appropriate distributive
structure seems important for general category theory,
especially in view of recent increased interest in braided
structures (as indicated, for example, by rough counting using
MathSciNet).
In this connection, our Referee~2 wrote this about the BD axioms:

\begin{quote}
``I \dots think that they probably are the `correct' notion of
braided distributivity (at least, for the intended applications)
but they are not at this stage actual coherence conditions.

Part of the reason for my caution on this matter is the fact
that these conditions are strongly derived from quantum-physical
considerations. \dots\  My concern is that there have been
previous situations where the natural physical generalization of
a concept and the natural categorical generalization of a
concept do not coincide."
\end{quote}

We leave the question of an appropriate coherence theorem open.

\begin{conj}\rm
There is a coherence theorem showing that the BD axioms are in
fact coherence conditions for braided distributivity.
\end{conj}

\begin{rmk}\rm
Our Referee~1 alerted us to the possibility of confusion between
our topic here and weak distributivity as developed in
\cite{CS}.  Both topics deal with connections between two
monoidal structures, and both have ``distributivity'' in their
names. They are nevertheless, in the words of the updated
version of \cite{CS}, ``virtually orthogonal''. As explained
there, it is very rare for distributive and weak distributive
laws to hold simultaneously. Also, in weak distributivity, the
two monoidal structures play interchangeable roles; in our
situation, just as in ordinary algebra, $\otimes$ distributes
over $\oplus$ but not vice versa.  Finally, our work is relevant
specifically for the situation where the monoidal structure
$\otimes$ is braided but not symmetric.
\end{rmk}

\section{BD categories} 
\label{sec:axioms}

The goal of this section is to present our notion of a category
equipped with a symmetric monoidal structure $\oplus$ and a
braided monoidal structure $\otimes$ with $\otimes$ distributing
over $\oplus$. We call such categories braided distributive, for
short BD, categories. The main task here is to present
appropriate requirements, or axioms, for the distributivity
isomorphisms.

Before giving a definition of BD categories, we
describe how our requirements were chosen, and we describe some
differences (both in content and in presentation) between our
work and that of Laplaza \cite{Laplaza}.

Given that we want certain isomorphisms (associativity, distributivity, etc.), how should we choose appropriate requirements for these isomorphisms to satisfy?  In the case of symmetric monoidal structures, as studied in \cite{MacLane1963} and \cite{Kelly}, the role of the requirements is to ensure the commutativity of all ``reasonable'' diagrams built from the associativity, commutativity, and unit isomorphisms and their inverses.  (``Reasonable'' requires careful formulation, for example to avoid expecting the special case $A\otimes A\to A\otimes A$ of commutativity to coincide with the identity morphism.)  In our situation, however, we do not want to equate the commutativity isomorphism $\gamma_{A,B}:A\otimes B\to B\otimes A$ with the inverse of $\gamma_{B,A}$, since that would make the braiding into a symmetry. So we have some freedom as to which diagrams should be required to commute; how can we responsibly exercise that freedom?

There are two mathematical constraints on this freedom, plus a
practical consideration that also influenced our choices.  The
first and most important mathematical constraint is that our
axioms should be satisfied in the examples that we set out to
describe, the modular tensor categories used in anyon models.
Axioms that fail in the intended examples are useless.  So we
must not make our axioms too strong.

The second mathematical constraint is that our axioms should not
be too weak; they should entail all the information needed in
our computations of specific examples. For instance, our axioms
should support the computations, as in \cite{G225}, of the
associativity and braiding matrices for Fibonacci anyons.

For practical purposes, we stay close to the requirements for
braided monoidal structures, as presented in \cite{JS1, JS2},
and for distributivity, as presented in \cite{Laplaza}. This
will enable us to make use of some of the computations in these
earlier papers.

\begin{rmk}\rm
For readers familiar with paper \cite{Laplaza} of Laplaza, we
point out two mathematical and three presentational differences
between his paper and what we do here.  The mathematical
differences are:
\begin{ls}
\item While Laplaza works with distributivity  given by
monomorphisms, our distributivity morphisms
$\delta_{A,B,C}:A\otimes(B\oplus C)\to (A\otimes
B)\oplus(A\otimes C)$ are assumed to be isomorphisms, not merely
monomorphisms.

  Recall that, according to Proviso~\ref{prov:laplaza}, speaking
  about Laplaza's requirements we mean, by default, the
  isomorphism-bound versions of his requirements.

\item On the other hand, Laplaza assumes that the multiplicative
monoidal structure is symmetric. We do not do that. Our
multiplicative monoidal structure, given by $\otimes$, is not
symmetric but only braided.

    In that sense our approach is more general then the
    isomorphism-bound case of Laplaza's approach.
\end{ls}
The presentational differences are these:
\begin{ls}
\item The dual distributivity morphisms, which Laplaza calls
$\delta^{\#}:(B\oplus C)\otimes A\to (B\otimes A)\oplus(C\otimes
A)$ are not taken as primitive data but are defined using
$\delta$ and the commutativity isomorphisms for $\otimes$.
\item Similarly, we do not take $0\otimes A\to0$ as primitive
but define it using $A\otimes0\to0$ and commutativity of
$\otimes$.
\item We take associativity isomorphisms in the
direction $(A\otimes B)\otimes C\to A\otimes(B\otimes C)$, and
similarly for $\oplus$, to agree with \cite{G225} and \cite{PP}.
Laplaza used the opposite direction.
\end{ls}
\end{rmk}

\begin{notat}\rm
We shall sometimes use the usual conventions from algebra that
$XY$ means $X\otimes Y$ and that, for example, $X\oplus YZ$
means $X\oplus(YZ)$, not $(X\oplus Y)Z$. For better fit with our
application \cite{G239}, we write composition of morphisms in
the left-to-right order, so $f \circ g$ means ``first $f$ and
then $g$''.
\end{notat}

After these preliminary comments, we now define BD categories.

\begin{df}\rm
A \emph{BD category} is a category \C\ equipped with
\begin{enumerate}
\item two bifunctors
$\oplus,\otimes: \C\times\C\to\C$ and two distinguished objects
$0,1$ in \C\ to serve as the units for the two bifunctors, and
\item isomorphisms
\begin{align*}
&\mathbf{associative}\ \oplus
&& \w{\alpha^{\oplus}_{A,B,C}}{(A\oplus B)\oplus C}{A\oplus(B\oplus C)}\\
&\mathbf{unit}\ \oplus
&&\w{\lambda^\oplus_A}{0\oplus A}{A}\quad\text{ and}\quad
   \w{\rho^\oplus_A}{A\oplus0}A\\
&\mathbf{commutative}\ \oplus
&&\w{\gamma^\oplus_{A,B}}{A\oplus B}{B\oplus A}\\
&\mathbf{associative}\ \otimes
&&\w{\alpha^{\otimes}_{A,B,C}}{(A\otimes B)\otimes C}{A\otimes(B\otimes C)}\\
&\mathbf{unit}\ \otimes
&&\w{\lambda^\otimes_A}{1\otimes A}A\quad\text{ and}\quad
\w{\rho^\otimes_A}{A\otimes1}A\\
&\mathbf{commutative}\ \otimes
&&\w{\gamma^\otimes_{A,B}}{A\otimes B}{B\otimes A}\\
&\mathbf{distributive}\ 2
&&\w{\delta_{A,B,C}}{A\otimes(B\oplus C)}{(A\otimes B)\oplus(A\otimes C)} \\
&\mathbf{distributive}\ 0
&&\w{\eps_A}{A\otimes 0}0
\end{align*}
natural with respect to the variable objects $A,B,C$, and subject to the axioms, the \emph{BD axioms}, given by Figures~1-18 below in this section.
\end{enumerate}
\end{df}

\begin{rmk}\rm\mbox{}
\begin{itemize}
\item 0 and 1 are not generally initial and terminal objects of \C, because $\oplus$ and $\otimes$ are not generally coproduct and product operations.
\item The 2 and 0 in the names of the distributivity isomorphisms refer to the number of summands on the right of $\otimes$; recall that the sum of no terms is understood to be 0.
\item If we worked with sets instead of categories and with equality instead of isomorphism, then these isomorphisms would provide the structure of a commutative semiring with unit. Here ``semi'' refers to the lack of additive inverses.
\end{itemize}
\end{rmk}

We next present the axioms in the form of commutative diagrams, but some words are needed about the form of these diagrams. Each of them is, when considered merely as an undirected graph, a simple cycle. But when the directions of the arrows are taken into account, they may not cohere, so it may not be immediately clear how to read such a diagram; hence the following explanation.  Any path $p$ in the undirected graph underlying such a diagram, say from vertex $X$ to vertex $Y$, represents the isomorphism from $X$ to $Y$ obtained by composing, in order along $p$, the isomorphisms written on the labels of the arrows that point in the same direction as $p$ and the inverses of the isomorphisms labeling arrows in the direction opposite to $p$. Less precisely but more memorably: As you go along $p$, compose the arrows you meet, but if the arrow points against the direction you're going then use the inverse isomorphism. Commutativity of a cycle diagram means that, for any two of its vertices $X$ and $Y$, both of the paths from $X$ to $Y$ represent the same isomorphism. It is easy to check that, if this happens for one choice of $X$ and $Y$, then it happens as well for all other choices of $X$ and $Y$. (It is permitted here that $X$ and $Y$ are the same vertex of the diagram; a path that goes all the way around the diagram must represent the identity isomorphism.)

We now present the desired axioms, in three groups: those that pertain to the additive structure, those that pertain to the multiplicative structure, and those that combine the two structures by distributivity.  In the captions under the diagrams, we name these requirements.

For the additive structure, we require the standard coherence conditions, in Figures~1 through 4, for a symmetric monoidal category, as given in \cite{Kelly}, simplifying an earlier version from \cite{MacLane1963}.

\begin{figure}[H]
\[\xymatrix@R+1pc{
& (A\oplus B)\oplus(C\oplus D)\ar@/^/[rd]^{\ap_{A,B,C\oplus D}}\\
((A\oplus B)\oplus C)\oplus D\ar@/^/[ru]^{\ap_{A\oplus B,C,D}}
           \ar[d]_{\ap_{A,B,C}\oplus1_D}
&& A\oplus(B\oplus(C\oplus D)) \\
(A\oplus(B\oplus C))\oplus D
    \ar[rr]_{\ap_{A,B\oplus C,D}}
&& A\oplus((B\oplus C)\oplus D) \ar[u]_{1_A\oplus\ap_{B,C,D}}
}\]
\caption{Additive Pentagon Condition}
\end{figure}

\begin{figure}[H]
\[\xymatrix@C+0pc@R+1pc{
&A\oplus(B\oplus C)\ar@/^/[r]^{\gp_{A,B\oplus C}}
&(B\oplus C)\oplus A\ar@/^/[rd]^{\ap_{B,C,A}}\\
(A\oplus B)\oplus C\ar@/^/[ru]^{\ap_{A,B,C}}
       \ar@/_/[rd]_{\gp_{A,B}\oplus 1_C}
&&&B\oplus(C\oplus A)\\
&(B\oplus A)\oplus C\ar@/_/[r]_{\ap_{B,A,C}}
&B\oplus(A\oplus C)\ar@/_/[ru]_{1_B\oplus\gp_{A,C}}
}\]
\caption{Additive Hexagon Condition}
\end{figure}

\begin{figure}[H]
\[\xymatrix@C+2pc@R+0pc{
(A\oplus0)\oplus B\ar@/^/[rr]^{\ap_{A,0,B}} \ar@/_/[rd]_{\rho^\oplus_A \oplus 1_B}
&&A\oplus(0\oplus B)\ar@/^/[ld]^{1_A\oplus\lambda^\oplus_B}\\
&A\oplus B
}\]
\caption{Additive Unit Associativity}
\end{figure}

\begin{figure}[H]
\[\xymatrix@C+8pc{
A\oplus B\ar@/^/[r]^{\gp_{A,B}} &B\oplus A\ar@/^/[l]^{\gp_{B,A}}
}\]
\caption{Additive Symmetry}
\end{figure}

For the multiplicative structure, we require the axioms in Figures~5 through 8 for a braided monoidal structure, as given by the coherence conditions of Joyal and Street in \cite{JS1,JS2}.  They are the same as for addition above except that symmetry is omitted and, to partially compensate for this omission, the hexagon condition is required to hold also when every $\gamma^\otimes_{X,Y}$ is replaced with ${\gamma^\otimes_{Y,X}}^{-1}$.

\begin{figure}[H]
\[\xymatrix@R+1pc{
& (A\x B)\x (C\x D)
  \ar@/^/[rd]^{\ax_{A,B,C\x D}}\\
((A\x B)\x C)\x D\ar@/^/[ru]^{\ax_{A\x B,C,D}}
           \ar[d]_{\ax_{A,B,C}\x 1_D}
&& A\x (B\x (C\x  D)) \\
(A\x (B\x C))\x D
    \ar[rr]_{\ax_{A,B\x C,D}}
&& A\x ((B\x C)\x D) \ar[u]_{1_A\x \ax_{B,C,D}}
}\]
\caption{Multiplicative Pentagon condition}
\end{figure}

\begin{figure}[H]
\[\xymatrix@C+0pc@R+1pc{
&A\x (B\x C)\ar@/^/[r]^{\gx_{A,B\x C}}
&(B\x C)\x A\ar@/^/[rd]^{\ax_{B,C,A}}\\
(A\x B)\x C\ar@/^/[ru]^{\ax_{A,B,C}}
       \ar@/_/[rd]_{\gx_{A,B}\x  1_C}
&&&B\x (C\x  A)\\
&(B\x A)\x C\ar@/_/[r]_{\ax_{B,A,C}}
&B\x (A\x C)\ar@/_/[ru]_{1_B\x \gx_{A,C}}
}\]
\caption{Multiplicative Hexagon:
  Moving one factor in front of two}
\end{figure}

\begin{figure}[H]
\[\xymatrix@C+0pc@R+1pc{
&A\x (B\x C)
&(B\x C)\x A\ar@/_/[l]_{\gx_{B\x C,A}}
  \ar@/^/[rd]^{\ax_{B,C,A}}\\
(A\x B)\x C\ar@/^/[ru]^{\ax_{A,B,C}}
&&&B\x (C\x  A)\ar@/^/[ld]^{1_B\x \gx_{C,A}}\\
&(B\x A)\x C \ar@/^/[lu]^{\gx_{B,A}\x  1_C}
  \ar@/_/[r]_{\ax_{B,A,C}}
&B\x (A\x C)
}\]
\caption{Multiplicative Hexagon:
  Moving one factor behind two}
\end{figure}

\begin{figure}[H]
\[\xymatrix@C+2pc@R+0pc{
(A\x1)\x B\ar@/^/[rr]^{\ax_{A,1\,B}}
  \ar@/_/[rd]_{\rho^\x_A \x 1_B}
&&A\x(1\x B)\ar@/^/[ld]^{1_A\x\lambda^\x_B}\\
&A\x B
}\]
\caption{Multiplicative Unit Associativity}
\end{figure}

In the names of the hexagon conditions, ``in front of'' and ``behind'' refer to the customary picture of braided commutativity in terms of geometric braids (the same picture that gave the name ``braided'' to
this weakening of symmetry).

Before turning to the last group of axioms, the ones that involve distributivity, we make some comments to relate our axioms to the isomorphism-bound version of Laplaza's coherence conditions in \cite{Laplaza}.  Laplaza's description of distributivity includes, along with morphisms for ``distributivity from the left'' (what we called $\delta_{A,B,C}$ and $\eps_A$ above), similar morphisms for distributivity from the right, $(B\oplus C)\otimes A\cong (B\otimes A)\oplus(C\otimes A)$ and $0\otimes A\cong0$.  His coherence conditions II and XV say that these right-distributivity morphisms are, as one might expect, obtainable from the left ones by means of commutativity of $\otimes$. We have chosen to take only the left morphisms as primitive and to regard the right ones as being defined from the left ones and commutativity.  (In effect, we have chosen to make the set of primitive isomorphisms small rather than symmetrical.)  So we no longer need Laplaza's coherence conditions II and XV.  But the fact that we are working with braided rather than symmetric multiplication leads to a complication here, namely, which version of commutativity shall we use to define the right distributivity isomorphisms in terms of the left ones?  Our first pair of axioms will say that it doesn't matter; both choices lead to the same right distributivity isomorphisms.  To express these requirements more clearly and succinctly, we introduce the following notation for the braiding.

\begin{notat}
  $\beta_{X,Y}=\gamma^\otimes_{X,Y}\circ\gamma^\otimes_{Y,X}$.
\end{notat}

Thus $\beta_{X,Y}$ is a isomorphism from $X\otimes Y$ to itself.  It would be just the identity $1_{X\otimes Y}$ if the multiplicative structure were symmetric; in general, it can be considered a measure of how far a braided structure deviates from symmetry.  Pictorially, if we imagine $\gamma^\otimes_{X,Y}$ as interchanging $X$ with $Y$ by moving $X$ in front of $Y$, then $\beta_{X,Y}$ moves $X$ all the way around $Y$ back to its initial position, first passing in front of $Y$ and then returning behind $Y$.

In terms of this $\beta$ notation, we can express our first axioms for distributivity as the pair of diagrams in Figure~9.

\begin{figure}[H]\small
\begin{minipage}{0.75\textwidth}
\[\xymatrix@C+4pc@R+0pc{
A\x(B\oplus C)\ar[r]^{\delta_{A,B,C}}
            \ar[d]^{\beta_{A,B\oplus C}}
&(A\x B)\oplus(A\x C)\ar[d]^{\beta_{A,B}\oplus\beta_{A,C}}\\
A\x(B\oplus C)\ar[r]^{\delta_{A,B,C}}
&(A\x B)\oplus(A\x C)}\]
\end{minipage}
\begin{minipage}{0.20\textwidth}
\[
\begin{tikzcd}
   A\x0\ar[out=130, in=50, loop, distance=6em]{}{\beta_{A,0}} &
\end{tikzcd}
\]
\end{minipage}
\caption{Right Distributive}
\end{figure}

The rest of our axioms for distributivity are given in
Figures~10 through 18.  They are essentially among Laplaza's
conditions, but rewritten in terms of our primitive
isomorphisms.
In the captions of these figures, we indicate the number of the corresponding condition in Laplaza's paper \cite{Laplaza}.

The first three of these requirement, Figures~10 through 12, say that distribution respects additive manipulations --- commutativity, associativity, and unit properties. That is, given $A\otimes S$ where $S$ is a sum, it doesn't matter whether we perform additive manipulations within $S$ and then apply distributivity or first apply distributivity and then perform the corresponding manipulations on the resulting sum.

\begin{figure}[H]
\[\xymatrix@C+8pc@R+1pc{
A(B\oplus C)\ar[r]^{\delta_{A,B,C}}
      \ar[d]_{1_A\x\gp_{B,C}}
&(AB)\oplus(AC)\\
A(C\oplus B)\ar[r]^{\delta_{A,C,B}}
&(AC)\oplus(AB)\ar[u]_{\gp_{AC,AB}}
}\]
\caption{Distribution Respects Additive Commutativity (Laplaza Cond.~I)}
\end{figure}

\begin{figure}[H]
\[\xymatrix@C+2pc@R+1pc{
A(B\oplus(C\oplus D))\ar@/^/[r]^{\delta_{A,B,C\oplus D}}
&AB\oplus A(C\oplus D)\ar@/^/[r]^{1_{AB}\oplus\delta_{A,C,D}}
&AB\oplus(AC\oplus AD)\\
A((B\oplus C)\oplus D) \ar@/_/[r]_{\delta_{A,B\oplus C,D}}
           \ar[u]^{{1_A}\x\ap_{B,C,D}}
& A(B\oplus C)\oplus AD\ar@/_/[r]_{\delta_{A,B,C}\oplus1_{AD}}
&(AB\oplus AC)\oplus AD\ar[u]_{\ap_{AB,AC,AD}}}
\]
\caption{Distribution Respects Additive Associativity (Laplaza Cond.~V)}
\end{figure}

\begin{figure}[H]
\[\xymatrix@C+8pc@R+1pc{
A(B\oplus0)     \ar[r]^{\delta_{A,B,0}}
           \ar[d]^{1_A\oplus\rho^\oplus_B}
&(AB)\oplus(A0) \ar[d]^{1_{AB}\oplus\eps_A}\\
AB
&(AB)\oplus0     \ar[l]^{\rho^\oplus_{AB}}
}\]
\caption{Distribution Respects 0 as Neutral (Laplaza Cond.~XXI)}
\end{figure}

Next are axioms saying that, when distributing a product of several factors across a sum, it doesn't matter whether one distributes the whole product at once or the individual factors one after the other.  The case of a product of two factors distributing across a sum of two summands is the obvious one; it implies (in the presence of the other axioms) the cases with more factors or summands. It is, however, also necessary to cover the cases where the number of factors or the number of summands is zero. So we get the four axioms in Figures~13 through 16. In our names for the axioms, the numbers 2 or 0 refer first to the number of factors and second to the number of summands.

\begin{figure}[H]
\[\xymatrix@C+6pc@R+1pc{
(AB)(C\oplus D)\ar[r]^{\ax_{A,B,C\oplus D}}
&A(B(C\oplus D))\ar[d]^{1_A\x\delta_{B,C,D}}\\
&A(BC\oplus BD)\ar[d]^{\delta_{A,BC,BD}}\\
(AB)C\oplus(AB)D \ar@{<-}[uu]^{\delta_{AB,C,D}}
&A(BC)\oplus A(BD)\ar@{<-}[l]^{\ax_{A,B,C}\oplus\ax_{A,B,D}}
}\]
\caption{Sequential Distribution $2\x2$ (Laplaza Cond.~VI)}
\end{figure}

\begin{figure}[H]
\[\xymatrix@C+8pc@R+1pc{
(AB)0  \ar[r]^{\ax_{A,B,0}}
       \ar[d]_{\eps_{AB}}
&A(B0) \ar[d]^{1_A\x\eps_B}\\
0
&A0    \ar[l]^{\eps_A}
}\]
\caption{Sequential Distribution $2\x0$ (Laplaza Cond.~XVIII)}
\end{figure}

\begin{figure}[H]
\[\xymatrix@C+2pc@R+1pc{
1(A\oplus B)      \ar@/^/[rr]^{\delta_{1,A,B}}
            \ar@/_/[rd]_{\lambda^\x_{A\oplus B}}
&&(1A)\oplus(1B) \ar@/^/[ld]^{\lambda^\x_A \oplus \lambda^\x_B}\\
&A\oplus B
}\]
\caption{Sequential Distribution $0\x2$ (Laplaza Cond.~XXIII)}
\end{figure}

\begin{figure}[H]
\[\xymatrix@C+9pc{
1\x0 \ar@/^/[r]^{\eps_1}
     \ar@/_/[r]_{\lambda^\x_0}
&0
}\]
\caption{Sequential Distribution $0\x0$ (Laplaza Cond.~XIV)}
\end{figure}

The remaining axioms for distributivity, in Figures~17 and 18, concern a product of two sums, like $(A\oplus B)(C\oplus D)$. Distributivity lets us expand this as a sum of four products, but there is a choice whether to apply distributivity first from the left, obtaining $((A\oplus B)C)\oplus((A\oplus B)D)$, or from the right, obtaining $(A(C\oplus D))\oplus(B(C\oplus D))$. One axiom (Figure~17) says that both choices produce the same final result, up to associativity and commutativity of addition. (Unfortunately, the associativity and commutativity make the diagram rather large, accounting for the lowest five arrows in Figure~17.  Furthermore, our decision to define right distributivity in terms of left distributivity plus commutativity enlarges the diagram even more. Specifically, in Figure~17, the path from the top vertex to the third vertex in the left column is just an instance of right distributivity, and the top three arrows in the right column amount to the sum of two more such instances.)

In addition, there are analogous but far simpler axioms for the case where one or both of the factors is the sum of no terms rather than of two.  Our labels for these requirements include numbers 2 or 0 indicating the number of summands in each factor.

\begin{figure}[H]\footnotesize
\[\xymatrix@C-3pc@R+0pc{
&(A\oplus B)(C\oplus D)\ar@/^/[rd]^{\delta_{A\oplus B,C,D}}
           \ar@/_/[ld]_{\gx_{A\oplus B,C\oplus D}}\\
(C\oplus D)(A\oplus B) \ar[d]_{\delta_{C\oplus D,A,B}}
&&(A\oplus B)C\oplus(A\oplus B)D\ar[d]^{\gx_{A\oplus B,C}\oplus\gx_{A\oplus B,D}}
\\
(C\oplus D)A\oplus(C\oplus D)B\ar@{<-}[d]_{\gx_{A,C\oplus D}\oplus\gx_{B,C\oplus D}}
&&C(A\oplus B)\oplus D(A\oplus B)\ar[d]^{\delta_{C,A,B}\oplus\delta_{D,A,B}}
\\
A(C\oplus D)\oplus B(C\oplus D)\ar[d]_{\delta_{A,C,D}\oplus\delta_{B,C,D}}
&&((CA)\oplus(CB))\oplus((DA)\oplus(DB)) \ar@{<-}[d]_%
  {(\gx_{A,C}\oplus\gx_{B,C})}^{\oplus\ (\gx_{A,D}\oplus\gx_{B,D})}
\\
((AC)\oplus(AD))\oplus((BC)\oplus(BD))\ar@{<-}[d]_{\ap_{AC\oplus AD,BC,BD}}
&&((AC)\oplus(BC))\oplus((AD)\oplus BD))
\\
(((AC)\oplus(AD))\oplus(BC))\oplus(BD)\ar@{<-}[d]_{\ap_{AC,AD,BC}\oplus1_{BD}}
&&(((AC)\oplus(BC))\oplus(AD))\oplus(BD)\ar[u]_{\ap_{AC\oplus BC,AD,BD}}
                        \ar[d]^{\ap_{AC,BC,AD}\oplus1_{BD}}
\\
((AC)\oplus((AD)\oplus(BC)))\oplus(BD)
&&((AC)\oplus((BC)\oplus(AD)))\oplus(BD)
  \ar@/^3pc/[ll]^{(1_{AC}\oplus\gp_{BC,AD})\oplus1_{BD}}
}\]
\caption{Expand $2\x2$ (Laplaza Cond.~IX)}
\end{figure}

\begin{figure}[H]\small
\begin{minipage}{0.75\textwidth}
\[\xymatrix@C+2pc@R+1pc{
(A\oplus B)0\ar[r]^{\gx_{A\oplus B,0}}
&0(A\oplus B)\ar[r]^{\delta_{0,A,B}}
&(0A)\oplus(0B)\ar@{<-}[d]^{\gx_{A,0}\oplus\gx_{B,0}}\\
0\ar@{<-}[u]^{\eps_{A\oplus B}}
&0\oplus0\ar[l]^{\lambda^\oplus_0}
&(A0)\oplus(B0)\ar[l]^{\eps_A \oplus \eps_B}}\]
\end{minipage}
\begin{minipage}{0.18\textwidth}
\[
\begin{tikzcd}
   0\x0\ar[out=130, in=50, loop, distance=6em]{}{\gx_{0,0}} &
\end{tikzcd}
\]
\end{minipage}
\caption{Expand $2\x0$ and $0\x0$ (Laplaza Conds.~XII and~X)}
\end{figure}

This completes our list of axioms.

\section{Deriving Laplaza's requirements} 
\label{sec:cons}

In this section, we verify that our distributivity requirements are as strong as (the isomorphism-bound version of) Laplaza's. Apart from the symmetric monoidal coherence conditions for $\oplus$ (Figures~1--4) and the braided monoidal coherence conditions for $\otimes$ (Figures~5--8), we have 12 requirements for distributivity; see Figures~9--18 and take into account that Figures~9 and 18 contain two requirements each.  But Laplaza has 24 coherence conditions. Nevertheless, we claim that our requirements imply all of Laplaza's, when these are suitably interpreted.

``Suitably interpreted'' here means simply that Laplaza's right distributivity morphisms $\delta^\#$ (for sums of two objects) and $\lambda^*$ (for the sum $0$ of no objects) are to be regarded not as primitive data but as defined from the left distributivity morphisms according to
\[
\delta^\#_{A,B,C}=\gamma^\otimes_{A\oplus B,C}\circ\delta_{C,A,B}\circ
(\gamma^\otimes_{A,C}\oplus\gamma^\otimes_{B,C})^{-1}
\]
and
\[
\lambda^*_A=\gamma^\otimes_{0,A}\circ\eps_A.
\]
We also record here the notational difference that Laplaza's $\rho^*$ is our $\eps$.

Some equivalent formulations of these definitions will be useful. First, in the definition of $\delta^\#$, the factor $(\gamma^\otimes_{A,C}\oplus\gamma^\otimes_{B,C})^{-1}$ can be replaced by ${\gamma^\otimes_{A,C}}^{-1}\oplus{\gamma^\otimes_{B,C}}^{-1}$; these are equal by functoriality of $\oplus$.

Second, the first two of our requirements, Right Distributivity (Figure~9), allow us to replace all of the $\gamma^\otimes$'s in the definitions of $\delta^\#$ and $\lambda^*$ with the inverses of $\gamma^\otimes$'s in which the two subscripts have been interchanged. Thus
\[
\delta^\#_{A,B,C}={\gamma^\otimes_{C,A\oplus B}}^{-1}\circ\delta_{C,A,B}\circ
(\gamma^\otimes_{C,A}\oplus\gamma^\otimes_{C,B})
\]
and
\[
\lambda^*_A={\gamma^\otimes_{A,0}}^{-1}\circ\eps_A.
\]

\begin{thm}\label{thm:strong}
The BD axioms imply all (isomorphism-bound versions of) Laplaza's coherence conditions.
\end{thm}

\begin{proof}
We need to verify that, with the interpretation of $\delta^\#$ and $\lambda^*$ (and $\rho^*$), described above, all 24 of Laplaza's coherence conditions are consequences of our axioms.  Our interpretation has already verified conditions II and XV from Laplaza's list; these just tell how $\delta^\#$ is related to $\delta$ and how $\lambda^*$ is related to $\rho^*$, via commutativity.

Another ten of Laplaza's conditions are included in our list, namely I (Fig.~10), V (Fig.~11), VI (Fig.~13), IX (Fig.~17), X and XII (Fig.~18), XIV (Fig.~16), XVIII (Fig.~14), XXI (Fig.~12), and XXIII (Fig.~15).

Laplaza provides \cite[pages~35--36]{Laplaza} nine implications between his coherence conditions and, after proving them, points out on pages~40--41 how they can be used to reduce the number of coherence conditions.  These results of Laplaza are, however, not enough to reduce all of his requirements to ours, and there are two reasons for this.

First, some of Laplaza's implications depend on the coherence assumptions for $\oplus$ and $\otimes$ individually. Although we have the same requirements for $\oplus$, we do not have the symmetry requirement for $\gamma^\otimes$.  We must therefore verify, in our setting with only braided $\otimes$, the implications that Laplaza obtains using symmetric $\otimes$. This affects the four items 3, 4, 5, and 9 in Laplaza's list of implications, so we must verify those four implications.

Second, Laplaza's list needs (the second item labeled 3 on page~41) at least two of requirements XVI, XVII, and XVIII; he has shown (item~7 on page~35, which uses only additive coherence) that any two of these three imply the remaining one.  We have assumed only one of the three, namely requirement~XVIII (our Sequential Distribution $2\otimes0$, Figure~14). So we should show that requirement~XVIII suffices to imply both requirements~XVI and XVII.  Of course, we need only show that XVIII implies one of XVI and XVII, because then Laplaza's proof gives us the remaining one.

So we have five implications to verify; we begin with the two easiest
ones.  Laplaza's proofs for  implications 5 and 9 use requirements for
$\otimes$ only in the form that, for any $X$, the three isomorphisms
$X1\cong X$ given by
\[
\rho^\otimes_X, \quad \gamma^\otimes_{X,1}\circ\lambda^\otimes_X,\quad
\text{and }\quad {\gamma^\otimes_{1,X}}^{-1}\circ\lambda^\otimes_X
\]
all coincide.  This instance of coherence does not need symmetry;
it is a consequence of the hexagon conditions, obtainable by replacing
two of the three objects in those conditions by 1 (see
\cite[Proposition~1]{JS1} or \cite[Proposition~2.1]{JS2}).

The remaining three implications, namely
\begin{align*}
  \text{VI}&\implies \text{VII},\\
\text{VI}&\implies \text{VIII}, \quad\text{and}\\
\text{XVIII}&\implies \text{XVII},
\end{align*}
require more work, as follows.  In each case, we shall first write
down the desired conclusion in the form of a diagram that we want to
commute. As before, if we ignore the directions of the arrows, these
diagrams will be simple cycles, and their commutativity is to be
understood as explained above in terms of the isomorphisms represented
by paths in the cycles.  We shall then gradually transform the desired
diagram into equivalent formulations until we reach a formulation that
is known to be true.

In preparation for the proof of $\text{VI}\implies \text{VII}$, we
record a property of braided monoidal categories proved by Joyal and
Street, \cite[Diagram~B5 on page~3]{JS1} and \cite[page~45]{JS2},
namely the commutativity of Figure~19.  It is a version of the
Yang-Baxter equation.

\begin{figure}[H]
\[\xymatrix@C+8pc@R+1pc{
(AB)C  \ar[r]^{\gx_{A,B}\x1_C}
&(BA)C \ar[r]^{\gx_{BA,C}}
&C(BA) \ar@{<-}[d]^{\ax_{C,B,A}}
\\
A(BC) \ar@{<-}[u]^{\ax_{A,B,C}}
&A(CB) \ar@{<-}[l]^{1_A\x\gx_{B,C}}
&(CB)A \ar@{<-}[l]^{\gx_{A,CB}}
}\]
\caption{Yang-Baxter Equation}
\end{figure}

We now turn to the proof that Laplaza's condition VI, which we called ``Sequential distribution $2\otimes2$'' (Figure~13), implies his condition VII. Recall that condition VI concerns two ways of applying distributivity to a product of the form $(AB)(C\oplus D)$ to obtain $(AB)C\oplus(AB)D$.  We can either apply distributivity directly, distributing $AB$ across $C\oplus D$, or we can, after using associativity, first distribute $B$ across $C\oplus D$ and then distribute $A$ across the resulting $BC\oplus BD$. Requirement~VI says that these two ways produce the same morphism.

Requirement~VII is the analog with multiplication in the other order.  That is, it asserts the equality of the morphism $(C\oplus D)(BA)\cong C(BA)\oplus D(BA)$ obtained by distributing $BA$ across the sum and the morphism obtained by first distributing $B$ and then $A$. Unfortunately, when written out in full, VII is considerably longer than VI, simply because the distributivity morphisms with the sum on the left, Laplaza's $\delta^\#$, are compositions of $\delta$'s and $\gamma$'s.  (Had we kept Laplaza's convention that $\delta^\#$ is primitive, VII would be shorter, but we would have to expand $\delta^\#$ in terms of $\delta$ during the proof, using II, and our work would be no easier.)  Figure 20 exhibits, in the solid arrows, the requirement~VII that we aim to prove.  The dashed arrows indicate where three of the solid arrows represent a single $\delta^{\#}$.

\begin{figure}[H]\small
\[\xymatrix@C+1pc@R+0pc{
(BA)(C\oplus D)       \ar@/^/[rr]^{\delta_{BA,C,D}}
&& (BA)C\oplus(BA)D  \ar[d]^{\gx_{BA,C}\oplus\gx_{BA,D}}
\\
(C\oplus D)(BA)       \ar@{<-}[u]^{\gx_{BA,C,D}}
                \ar@{-->}[rr]^{\delta^\#_{C,D,BA}}
&&C(BA)\oplus D(BA)   \ar@{<-}[d]^{\ax_{C,B,A}\oplus\ax_{D,B,A}}
\\
((C\oplus D)B)A       \ar[u]^{\ax_{C\oplus D,B,A}}
                \ar@{-->}[rddd]^{\delta^\#_{C,D,B}\x1_A}
&&(CB)A\oplus(DB)A   \ar@{<-}[d]^{\gx_{A,CB}\oplus\gx_{A,DB}}
\\
(B(C\oplus D))A       \ar[u]^{\gx_{B,C\oplus D}\x1_A}
&&A(CB)\oplus A(DB)   \ar@{<-}[d]^{\delta_{A,CB,DB}}
\\
(BC\oplus BD)A    \ar@{<-}[u]^{\delta_{B,C,D}\x1_A}
&&A(CB\oplus DB)  \ar@/^/[ld]^-{\gx_{A,CB\oplus DB}}
\\
&(CB\oplus DB)A   \ar@/^/@{<-}[lu]^-{(\gx_{B,C}\oplus\gx_{B,D})\x1_A}
                \ar@{-->}[ruuu]^{\delta^\#_{CB,DB,A}}
}\]
\caption{Laplaza Cond.~VII}
\end{figure}

We now begin to transform Figure~20, using Figure~13 and the
requirements on $\otimes$. Consider first the top
three (solid) arrows in the right column of Figure~20.  Each of these
involves a sum of two morphisms, one working with $C$ and one with
$D$.  Because of the functoriality of $\oplus$, we can treat these three
arrows for each summand separately. The three that work with $C$
represent  the composite morphism
\[
\gamma^\otimes_{BA,C}\circ{\alpha^\otimes_{C,B,A}}^{-1}\circ
{\gamma^\otimes_{A,CB}}^{-1},
\]
which is exactly the morphism represented in the Yang-Baxter
equation by the path from $(BA)C$ to $A(CB)$ that goes around the
right side of the diagram in Figure~19.  So we can rewrite this
morphism as given by the path from $(BA)C$ to $A(CB)$ that goes
around the left side of Figure~19, namely
\[
({\gamma^\otimes_{A,B}}^{-1}\otimes 1_C)\circ \alpha^\otimes_{A,B,C}\circ
(1_A\otimes\gamma^\otimes_{B,C}).
\]
The same replacement can be done for the other summand in the top
three arrows on the right of Figure~20, since the same argument
applies with $D$ in place of $C$.

Before recording the result of this transformation, we observe another
transformation that can be performed on a disjoint part of Figure~20. In the
lower left part, we have three morphisms of the form
$\xi\otimes1_A$, where $\xi$ involves $B$, $C$, and $D$. Preceding
those three morphisms in clockwise order is a $\gamma^\otimes$ that
moves $A$ from the left to the right without affecting the part that
involves $B$, $C$, and $D$.  By naturality of $\gamma$, we can,
without changing the composition of these morphisms, move $A$ to the
other side of the $B$, $C$, and $D$ part after the other three
morphisms rather than before.

The equivalent form of VII obtained by these transformations is shown in Figure~21.

\begin{figure}[H]\small
\[\xymatrix@C+1pc@R+0pc{
(BA)(C\oplus D)       \ar@/^/[rr]^{\delta_{BA,C,D}}
&& (BA)C\oplus(BA)D  \ar@{<-}[d]^{\gx_{A,B}\x1_C\oplus\gx_{A,B}\x1_D}
\\
(C\oplus D)(BA)       \ar@{<-}[u]^{\gx_{BA,C,D}}
&&(AB)C\oplus(AB)D   \ar[d]^{\ax_{A,B,C}\oplus\ax_{A,B,D}}
\\
((C\oplus D)B)A       \ar[u]^{\ax_{C\oplus D,B,A}}
&&A(BC)\oplus A(BD)   \ar[d]^{1_A\x\gx_{B,C}\oplus1_A\x\gx_{B,D}}
\\
A((C\oplus D)B)      \ar[u]^{\gx_{A,(C\oplus D)B}}
&&A(CB) \oplus A(DB)\ar@{<-}[d]^{\delta_{A,CB,DB}}
\\
A(B(C\oplus D))      \ar[u]^{1_A\x\gx_{B,C\oplus D}}
&&A(CB\oplus DB) \ar@{<-}@/^/[ld]^-{\quad1_A\x(\gx_{B,C}\oplus\gx_{BD})}
\\
&A(BC\oplus BD)  \ar@/^/@{<-}[lu]^-{1_A\x\delta_{B,C,D}\ }
}\]
\caption{\IS{1}{VII}}
\end{figure}

Here the three morphisms in the lower right can, by naturality of
$\delta$, be equivalently replaced with a single arrow labeled
$\delta_{A,BC,BD}$, pointing from the $A((BC)\oplus(BD))$ at the bottom to
the $A(BC)\oplus A(BD)$ at the middle of the right column.  This arrow
and the two adjacent arrows in Figure~21 match the lower path in Figure~13
from the lower left to the upper right corner (clockwise in Figure~21
matches counterclockwise in Figure~13).  So these arrows in Figure~21
can be equivalently replaced by the other path between the same
corners of Figure~13.  The result is Figure~22.

\begin{figure}[H]\small
\[\xymatrix@C+6pc@R+0pc{
(BA)(C\oplus D)       \ar@/^/[r]^{\delta_{BA,C,D}}
& (BA)C\oplus(BA)D  \ar@{<-}[d]^{\gx_{A,B}\x1_C\oplus\gx_{A,B}\x1_D}
\\
(C\oplus D)(BA)       \ar@{<-}[u]^{\gx_{BA,C,D}}
&(AB)C\oplus(AB)D   \ar@{<-}[d]^{\delta_{AB,C,D}}
\\
((C\oplus D)B)A       \ar[u]^{\ax_{C\oplus D,B,A}}
&(AB)(C\oplus D)     \ar[ldd]^{\ax_{A,B,C\oplus D}}
\\
A((C\oplus D)B)      \ar[u]^{\gx_{A,(C\oplus D)B}}
\\
A(B(C\oplus D))      \ar[u]^{1_A\x\gx_{B,C\oplus D}}
}\]
\caption{\IS{2}{VII}}
\end{figure}

The arrow across the top of Figure~22 and the two arrows in the right
column can be equivalently replaced by a single arrow, thanks to the
naturality of $\delta$ (applied to the morphisms $\gamma_{A,B}$,
$1_C$, and $1_D$).  The result of this simplification is Figure~23.

\begin{figure}[H]\small
\[\xymatrix@C+6pc@R+0pc{
(BA)(C\oplus D)       \ar@{<-}[rdd]^{\gx_{A,B}\x1_{C\oplus D}}
\\
(C\oplus D)(BA)       \ar@{<-}[u]^{\gx_{BA,C,D}}
\\
((C\oplus D)B)A       \ar[u]^{\ax_{C\oplus D,B,A}}
&(AB)(C\oplus D)     \ar[ldd]^{\ax_{A,B,C\oplus D}}
\\
A((C\oplus D)B)      \ar[u]^{\gx_{A,(C\oplus D)B}}
\\
A(B(C\oplus D))      \ar[u]^{1_A\x\gx_{B,C\oplus D}}
}\]
\caption{\IS{3}{VII}}
\end{figure}

This diagram, which no longer involves $\delta$, commutes by virtue of
the Yang-Baxter equation (Figure~19).  This completes the proof that
Laplaza's requirement~VII follows from VI plus the coherence
properties of braided monoidal categories plus naturality of $\delta$.
\end{proof}

We turn next to the proof that $\text{VI}\implies\text{VIII}$. Laplaza's requirement~VIII, our goal in this proof, is, like VI and VII, about distributing two factors, $A$ and $B$, across a sum $C\oplus D$, but now one factor $A$ is on the left of $C\oplus D$ and the other factor $B$ is on the right. The requirement says that we get the same isomorphism whether we distribute $A$ across $C\oplus D$ first and then $B$ across $(AC)\oplus(AD)$ or we distribute $B$ across $C\oplus D$ first and then $A$ across $(CB)\oplus(DB)$.  Figure~24 is the diagram whose commutativity expresses this.  As before, we have indicated with dashed arrows the places where three of our arrows represent a right distributivity isomorphism $\delta^{\#}$.

\begin{figure}[H]\small
\[\xymatrix@C+1pc@R+0pc{
A(BC\oplus BD)       \ar@/^2pc/[r]^{1_A\x(\gx_{B,C}\oplus\gx_{B,D})}
&A(CB\oplus DB)      \ar@/^2pc/[r]^{\delta_{A,CB,DB}}
&A(CB)\oplus A(DB)      \ar@{<-}[d]^{\ax_{A,C,B}\oplus\ax_{A,D,B}}
\\
A(B(C\oplus D))          \ar[u]^{1_A\x\delta_{B,C,D}}
&&(AC)B\oplus(AD)B      \ar@{<-}[d]^{\gx_{B,AC}\oplus\gx_{B,AD}}
\\
A((C\oplus D)B)          \ar@{<-}[u]^{1_A\x\gx_{B,C\oplus D}}
                   \ar@{-->}[ruu]_{1_A\x\delta^\#_{C,D,B}}
&&B(AC)\oplus B(AD)      \ar@{<-}[d]^{\delta_{B,AC,AD}}
\\
(A(C\oplus D))B         \ar[u]^{\ax_{A,C\oplus D,B}}
&(AC\oplus AD)B     \ar@{<-}@/^2pc/[l]^{\delta_{A,C,D}\x1_B}
                  \ar@{-->}[ruu]^{\delta^\#_{AC,AD,B}}
&B(AC\oplus AD)     \ar@/^2pc/[l]^{\gx_{A,(AC)\oplus(AD)}}
}\]
\caption{Laplaza Cond.~VIII}
\end{figure}

We now proceed, as in the preceding proof, to transform this diagram
into equivalent diagrams, until we obtain one that is known to
commute.  We begin with the two arrows across the top of the diagram,
which can, by naturality of $\delta$, be equivalently replaced by
\[
\delta_{A,(BC)\oplus(BD)}\circ
\Big((1_A\otimes\gamma^\otimes_{B,C})\oplus(1_A\otimes\gamma^\otimes_{B,D})\Big).
\]
The first factor here, $\delta_{A,(BC)\oplus(BD)}$, together with the top
arrow in the left column of Figure~24, constitutes a path in the
Sequential distribution 22 condition, Figure~13. So these two arrows
can be replaced with the other path between the same vertices of
Figure~13. The resulting diagram, equivalent to Figure~24, is
Figure~25.

\begin{figure}[H]\small
\[\xymatrix@C+0pc@R+0pc{
(AB)C\oplus(AB)D         \ar[rr]^{\ax_{A,B,C}\oplus\ax_{A,B,D}}
&&A(BC)\oplus A(BD)       \ar[d]^{(1_A\x\gx_{B,C})\oplus(1_A\x\gx_{B,D})}
\\
(AB)(C\oplus D)         \ar[u]^{\delta_{AB,C,C}}
&&A(CB)\oplus A(DB)      \ar@{<-}[d]^{\ax_{A,C,B}\oplus\ax_{A,D,B}}
\\
A(B(C\oplus D))          \ar@{<-}[u]^{\ax_{A,B,C\oplus D}}
&&(AC)B\oplus(AD)B      \ar@{<-}[d]^{\gx_{B,AC}\oplus\gx_{B,AD}}
\\
A((C\oplus D)B)          \ar@{<-}[u]^{1_A\x\gx_{B,C\oplus D}}
&&B(AC)\oplus B(AD)      \ar@{<-}[d]^{\delta_{B,AC,AD}}
\\
(A(C\oplus D))B         \ar[u]^{\ax_{A,C\oplus D,B}}
&(AC\oplus AD)B     \ar@{<-}@/^2pc/[l]^{\delta_{A,C,D}\x1_B}
&B(AC\oplus AD)     \ar@/^2pc/[l]^{\gx_{A,(AC)\oplus(AD)}}
}\]
\caption{\IS{1}{VIII}}
\end{figure}

Next, consider the arrow across the top and the upper three arrows in
the right column of Figure~25. Each is the sum of two
isomorphisms; temporarily concentrate on the first summand in each,
i.e., the part involving $C$ rather than $D$.  So we are looking at
the composite
\[
\alpha^\otimes_{A,B,C} \circ
(1_A\otimes\gamma^\otimes_{B,C})\circ
{\alpha^\otimes_{A,C,B}}^{-1}\circ
{\gamma^\otimes_{B,AC}}^{-1}.
\]
This is, up to interchanging the roles of $A$ and $B$, represented by
a path of four arrows in the first Multiplicative Hexagon Condition,
Figure~6. So it can be replaced by the other path joining the same end
vertices, a path of length two.  This simplification of the summands
involving $C$ applies equally well to the summands involving $D$. By
functoriality of $\oplus$, we can combine the simplifications to put our
Figure~25 into the equivalent form Figure~26.

\begin{figure}[H]\small
\[\xymatrix@C+0pc@R+0pc{
(AB)C\oplus(AB)D         \ar@{<-}[rrd]^{\qquad(\gx_{B,A}\x1_C)\oplus(\gx_{B,A}\x1_D)}
\\
(AB)(C\oplus D)           \ar[u]^{\delta_{AB,C,C}}
&&(BA)C\oplus(BA)D       \ar[dd]^{\ax_{B,A,C}\oplus\ax_{B,A,D}}
\\
A(B(C\oplus D))           \ar@{<-}[u]^{\ax_{A,B,C\oplus D}}
\\
A((C\oplus D)B)           \ar@{<-}[u]^{1_A\x\gx_{B,C\oplus D}}
&&B(AC)\oplus B(AD)       \ar@{<-}[d]^{\delta_{B,AC,AD}}
\\
(A(C\oplus D))B          \ar[u]^{\ax_{A,C\oplus D,B}}
&(AC\oplus AD)B      \ar@{<-}@/^2pc/[l]^{\delta_{A,C,D}\x1_B}
&B(AC\oplus AD)      \ar@/^2pc/[l]^{\gx_{A,(AC)\oplus(AD)}}
}\]
\caption{\IS{2}{VIII}}
\end{figure}

Now the two morphisms across the bottom of this figure can be replaced, thanks to naturality of $\gamma^\otimes$, by $\gamma^\otimes_{B,A(C\oplus D)}$ on the left and $1_B\otimes{\delta_{A,C\oplus D}}^{-1}$ on the right. The former and the bottom three arrows in the left column constitute a path of length four in the first Multiplicative Hexagon Condition, Figure~6 (with $B$, $A$, and $C\oplus D$ in the roles of $A$, $B$, and $C$, respectively), so they can be equivalently replaced by the other path between the same vertices. The result is Figure~27.

\begin{figure}[H]\small
\[\xymatrix@C+2pc@R+0pc{
(AB)C\oplus(AB)D         \ar@{<-}[rrd]^{\qquad(\gx_{B,A}\x1_C)\oplus(\gx_{B,A}\x1_D)}
\\
(AB)(C\oplus D)           \ar[u]^{\delta_{AB,C,C}}
&&(BA)C\oplus(BA)D       \ar[d]^{\ax_{B,A,C}\oplus\ax_{B,A,D}}
\\
(BA)(C\oplus D)           \ar[u]^{\gx_{B,A}\x1_{C\oplus D}}
&&B(AC)\oplus B(AD)       \ar@{<-}[d]^{\delta_{B,AC,AD}}
\\
B(A(C\oplus D))           \ar@{<-}[u]^{\ax_{B,A,C\oplus D}}
&&B(AC\oplus AD)      \ar@{<-}[ll]^{1_B\x\delta_{A,C,D}}
}\]
\caption{\IS{3}{VIII}}
\end{figure}

By naturality of $\delta$ and functoriality of $\oplus$ (specifically, $1_{C\oplus D}=1_C\oplus1_D$), we can replace the top arrow here and the upper two arrrows in the left column by a single arrow, labeled $\delta_{BA,C,D}$, pointing from the $(BA)(C\oplus D)$ in the left column to the $(BA)C\oplus(BA)D$ at the top of the right column.

But the resulting diagram is Figure~13, so the proof of
$\text{VI}\implies\text{VIII}$ is complete.

It remains to prove the implication
$\text{XVIII}\implies\text{XVII}$ which goes beyond Laplaza's redundancies.

\begin{proposition}\label{redundancy}
Laplaza's requirement XVIII implies his requirement XVII both in his original context and in the isomorphism-bound version of his context.
\end{proposition}

\begin{proof}[Proof of the proposition]
Actually, the two contexts coincide for the present purpose, because Laplaza requires the relevant distributivity maps $\lambda^*$ and $\rho^*$ to be isomorphisms.

We note that XVIII and XVII are the analogs of VI and VIII, respectively, where the sum of two terms $C\oplus D$ in VI and VIII is replaced by the sum of no terms, 0, in XVIII and XVII. Our proof of $\text{VI}\implies\text{VIII}$ can be made into a proof of $\text{XVIII}\implies\text{XVII}$ by systematically replacing binary sums with nullary sums throughout the computation. Some parts of the computation simplify, and we present here the resulting proof of $\text{XVIII}\implies\text{XVII}$.

In our notation, the goal XVII is the commutativity of
Figure~28. (Recall, in this connection, that Laplaza's $\rho^*_A$ is
our $\eps_A$ and that Laplaza's $\lambda^*_A$ is our
$\gamma^\otimes_{0,A}\circ \eps_A$ and also our
${\gamma^\otimes_{A,0}}^{-1}\circ \eps_A$, the latter two being equal by
requirement Right Distributive 0.)

\begin{figure}[H]
\[\xymatrix@C+3pc@R+0pc{
A(0B)         \ar@{<-}[d]_{1_A\x\gx_{A,0}}
&&(A0)B       \ar[ll]_{\ax_{A,0,B}}
\\
A(B0)         \ar[d]_{1_A\x\eps_B}
&&0B           \ar@{<-}[u]_{\eps_A\x1_B}
\\
A0            \ar@/_1pc/[r]_{\eps_A}
&0            \ar@/_1pc/@{<-}[r]_{\eps_B}
&B0           \ar[u]_{\gx_{B,0}}
}\]
\caption{Laplaza Cond.~XVII}
\end{figure}

Using naturality of $\gamma^\otimes$, we can replace the right column
as shown in Figure~29.

\begin{figure}[H]
\[\xymatrix@C+3pc@R+0pc{
A(0B)         \ar@{<-}[d]_{1_A\x\gx_{A,0}}
&&(A0)B       \ar[ll]_{\ax_{A,0,B}}
\\
A(B0)         \ar[d]_{1_A\x\eps_B}
&&B(A0)       \ar[u]_{\gx_{B,A0}}
\\
A0            \ar@/_1pc/[r]_{\eps_A}
&0            \ar@/_1pc/@{<-}[r]_{\eps_B}
&B0           \ar@{<-}[u]_{1_B\x\eps_{A}}
}\]
\caption{\IS{1}{XVII}}
\end{figure}

The top half of this diagram, i.e., the arrow across the top and the
upper arrows in both the left and the right columns, form a path in an
instance of the first Multiplicative Hexagon Condition, Figure~6. So we can
equivalently replace them by the other path joining the same vertices,
obtaining Figure~30.

\begin{figure}[H]
\[\xymatrix@C+3pc@R+0pc{
(AB)0         \ar[d]_{\ax_{A,B,0}}
&&(BA)0       \ar[ll]_{\gx_{B,A}\x1_0}
\\
A(B0)         \ar[d]_{1_A\x\eps_B}
&&B(A0)       \ar@{<-}[u]_{\ax_{B,A,0}}
\\
A0            \ar@/_1pc/[r]_{\eps_A}
&0            \ar@/_1pc/@{<-}[r]_{\eps_B}
&B0           \ar@{<-}[u]_{1_B\x\eps_{A}}
}\]
\caption{\IS{2}{XVII}}
\end{figure}

Thanks to Figure~14, the path on the left from $(AB)0$ to $0$ represents
$\eps_{AB}$ and the similar path on the right represents
$\eps_{BA}$. So our diagram simplifies to Figure~31.

\begin{figure}[H]
\[\xymatrix@C+3pc@R+0pc{
(AB)0         \ar[rd]_{\eps_{AB}}
&&(BA)0       \ar[ll]_{\gx_{B,A}\x1_0}
\\
&0            \ar@{<-}[ru]_{\eps_{BA}}
}\]
\caption{\IS{3}{XVII}}
\end{figure}

This diagram commutes by naturality of $\eps$, and so the proof of
$\text{XVIII}\implies\text{XVII}$ is complete.

This concludes the proof of the theorem.
\end{proof}

\section{Abelian categories plus distributivity} 
\label{sec:abel}
\renewcommand{\thefigure}{\hspace{-.333333em}}

By Theorem~\ref{thm:strong}, the BD axioms are strong enough to imply suitable versions of Laplaza's conditions from \cite{Laplaza}. In the present section, our goal is to show that the BD axioms are not excessively strong: They hold in the intended applications, the categories used to model non-abelian anyons. In fact, the BD axioms are consequences of just a part of the axiom system for modular tensor categories as described in \cite{PP} and \cite{G225}. Specifically, we shall use the following assumptions about a category \C\ to deduce the framework described in Section~\ref{sec:axioms}.

\begin{enumerate}
\item[A1.] \scr C is an abelian category.
\item[A2.] \scr C has a braided monoidal structure, with product operation
  $\otimes$.
\item[A3.] The bifunctor $\otimes$ is additive in each of its two arguments.
\end{enumerate}

We begin by summarizing some background information that we shall need
in our proofs. We use Chapter II of Freyd's book \cite{Freyd} as a
reference for the facts we need about abelian categories. The most
important of these facts, for our purposes, are the following.  An
abelian category has all finite products and coproducts, and these
coincide, i.e., there is a \emph{zero} object $0$ that is
simultaneously terminal and initial, and any two objects $A$ and $B$
have a \emph{sum} $A\oplus B$ that is simultaneously their product
with projections $p_A:A\oplus B\to A$ and $p_B:A\oplus B\to B$ and
their coproduct with injections $u_A:A\to A\oplus B$ and $u_B:B\to
A\oplus B$. For each two objects $A$ and $B$, the set
$\Hom_\C(A,B)$ of morphisms from $A$ to $B$ has the
structure of an abelian group; its group operation will be written as
$+$, and its zero element $0$ is the unique morphism from $A$ to $B$
that factors through the zero object. The composition operation of
\scr C is additive in each of its two arguments.

In addition to these consequences of assumption~A1, the availability of finite products in \scr C provides a symmetric monoidal structure; see Sections~VII.1 and VII.7 of \cite{MacLane1971}.  That is, there exist natural transformations $\alpha^\oplus$, $\lambda^\oplus$, $\rho^\oplus$, and $\gamma^\oplus$ satisfying the
BD axioms in Figures~1 through 4.

For future reference, we record here how these natural transformations are defined, in terms of the products in \scr C.  Two of them are trivial,\\ since $\lambda^\oplus_A : 0\oplus A\to A$ and $\rho^\oplus_A : A\oplus0\to A$ are simply the projections of the products to the second and first factor, respectively. (The inverses are defined as the unique morphisms into the products that act as $1_A$ to one factor and the unique morphism to the terminal object to the other factor.)

To describe $\alpha^\oplus$ and $\gamma^\oplus$, it will be useful to have a
very precise notation for the projections of a (binary) product; the notation
will also be useful in subsequent calculations.
We shall write $p_{A,B,1}$ for
the projection $A\oplus B\to A$ to the first factor and $p_{A,B,2}$ for the
projection $A\oplus B\to B$ to the second factor. If $A$ and $B$ are
sufficiently clear from the context, we may write simply $p_1$ and $p_2$.
Alternatively, we may use the abbreviated notations $p_A$ and $p_B$, as long as
$A$ and $B$ are distinct. But in general, the full three-subscript notation
serves to eliminate any danger of ambiguity.

With this notation, we can describe the associativity and
commutativity isomorphisms by telling how they compose with
projections. This information will suffice to completely determine
those isomorphisms, because a morphism into a product is determined by
its composites with projections.  For commutativity
$\gamma^\oplus_{A,B}:A\oplus B\to B\oplus A$, we have
\begin{align*}
\gamma^\oplus_{A,B}\circ p_{B,A,1}\ &=\ p_{A,B,2}:A\oplus B\to B,\\
\gamma^\oplus_{A,B}\circ p_{B,A,2}\ &=\ p_{A,B,1}:A\oplus B\to A. \end{align*}
For associativity $\alpha^\oplus_{A,B,C}:(A\oplus B)\oplus C\to
A\oplus(B\oplus C)$, we have
\begin{align*}
&\alpha^\oplus_{A,B,C}\circ p_{A,B\oplus C,1}&&=\
p_{A\oplus B,C,1}\circ p_{A,B,1}&:&&
(A\oplus B)\oplus C\to A\\
&\alpha^\oplus_{A,B,C}\circ p_{A,B\oplus C,2}\circ p_{B,C,1}&&=\
p_{A\oplus B,C,1}\circ p_{A,B,2}&:&&
(A\oplus B)\oplus C\to B\\
&\alpha^\oplus_{A,B,C}\circ p_{A,B\oplus C,2}\circ p_{B,C,2}&&=\
p_{A\oplus B,C,2}&:&&
(A\oplus B)\oplus C\to C
\end{align*}
The inverse isomorphisms admit similar descriptions.

So far, we have used only assumption~A1, that \C\ is an abelian category.  Assumption~A2 provides the multiplicative structure, with multiplication $\otimes$, unit 1, and associativity, commutativity, and unit isomorphisms satisfying the BD axioms in Figures~5 through 8.

To obtain distributivity isomorphisms and their axioms,
we use assumption~A3 as follows.  The projections $p_i:X_1\oplus X_2\to
X_i$ and injections $u_i:X_i\to X_1\oplus X_2$ of a sum satisfy the
equations
\begin{align*}
&u_i\circ p_i &&=\ 1_{X_i}  \\
&u_i\circ p_j &&=\ 0\phantom{mmmmmm}\text{for }i\neq j\\
&(p_1\circ u_1)+(p_2\circ u_2) &&=\ 1_{X_1\oplus X_2}
\end{align*}
where the addition is the group operation in $\Hom_\C(X_1\oplus X_2,X_1\oplus X_2)$.  Furthermore, these equations characterize sums, in the sense that any object $Y$ equipped with morphisms $p_i:Y\to X_i$ and $u_i:X_i\to Y$ satisfying these equations is canonically isomorphic to $X_1\oplus X_2$ (with the $p$'s and $u$'s for $Y$ corresponding, via the isomorphism, to those for $X_1\oplus X_2$). See Theorem~2.41 of \cite{Freyd}.

Consider now a sum $B\oplus C$, with its two projections $p_{B,C,i}$ and injections $u_{B,C,i}$ (in an obvious notation).  For any object $A$, the functor $A\otimes -$ is additive, by A3. Therefore the morphisms $A\otimes p_{B,C,i}$ and $A\otimes u_{B,C,i}$ also satisfy the sum equations and thus make $A\otimes(B\oplus C)$ canonically isomorphic to $(A\otimes B)\oplus(A\otimes C)$.  This canonical isomorphism will serve as our distributivity isomorphism $\delta_{A,B,C}$.  Writing out in detail the definition of this canonical isomorphism, we have
\begin{align*}
\delta_{A,B,C}\circ p_{A\otimes B,A\otimes C,1}
&\ =\ 1_A\otimes p_{B,C,1}: A\otimes(B\oplus C)\to A\otimes B,\\
\delta_{A,B,C}\circ p_{A\otimes B,A\otimes C,2}
&\ =\ 1_A\otimes p_{B,C,2}: A\otimes(B\oplus C)\to A\otimes C.
\end{align*}

We should check that $\delta$ is a natural transformation, since this
is one of the requirements of BD categories. So we must
check that, for any morphisms $a:A\to A'$, $b:B\to B'$, and $c:C\to
C'$, the composites
\[
\delta_{A,B,C}\circ ((a\otimes b)\oplus(a\otimes c))
\]
and
\[
(a\otimes(b\oplus c))\circ\delta_{A',B',C'}
\]
coincide. Since these are morphisms into a sum $(A'\otimes B')\oplus(A'\otimes C')$, it suffices to check that their composites with the projections to $A'\otimes B'$ and $A'\otimes C'$ coincide. We check the first of these; the second is entirely analogous. We compute for the first of the two allegedly coinciding morphisms
\begin{align*}
&\delta_{A,B,C}\circ ((a\otimes b)\oplus(a\otimes c))\circ
p_{A'\otimes B',A'\otimes C',1}=\\
&\delta_{A,B,C}\circ p_{A\otimes B,A\otimes C,1}\circ(a\otimes b)=\\
&(1_A\otimes p_{B,C,1})\circ(a\otimes b)
\end{align*}
and for the second
\begin{align*}
&(a\otimes(b\oplus c))\circ\delta_{A',B',C'}\circ
p_{A'\otimes B',A'\otimes C',1}=  \\
&(a\otimes(b\oplus c))\circ (1_{A'}\otimes p_{B',C',1})=\\
&a\otimes(p_{B,C,1}\circ b).
\end{align*}
(In the first of these computations, we used the definition of how the product bifunctor $\oplus$ acts on morphisms, and then we used the definition of $\delta$. In the second computation, we used the definition of $\delta$ and then the fact that $\otimes$ is functorial in both arguments along with the definition of $\oplus$ on morphisms.) The last lines in our two computations agree because $\otimes$ is a bifunctor.  This completes the verification that $\delta$ is natural.

We still need to specify the nullary part of the distributivity structure, the natural isomorphisms $\eps_A:A\otimes0\to0$, but this is trivial; 0 is terminal so any object has just one morphism to 0. (We could equivalently say $\eps_A=p_{A,0,2}$.) Naturality of $\eps$ is trivial also, since it requires equality between morphisms to the terminal object.  To show that $\eps_A$ is an isomorphism, we consider the unique morphism $0\to A\otimes0$ given by initiality of $0$. Composing it with $\eps_A$ in one order, we get an endomorphism of 0, which is the identity of 0 because 0 is terminal (and also because 0 is initial). The composite in the other order is an endomorphism of $A\otimes0$ that factors through 0 and is therefore the zero element of the endomorphism group $Hom(A\otimes0,A\otimes0)$. We must show that this zero endomorphism is also the identity endomorphism. But this property, ``zero equals identity'' is true for the object 0 and is preserved by $A\otimes-$ because this is functorial and bilinear.

Having specified all the natural isomorphisms required in the definition of BD categories, we must still verify all the BD axioms.

\begin{thm}
Let \C\ be a category satisfying the conditions A1--A3 above, and let natural isomorphisms
\begin{align*}
& \w{\alpha^{\oplus}_{A,B,C}}{(A\oplus B)\oplus C}{A\oplus(B\oplus C)}\\
&\w{\lambda^\oplus_A}{0\oplus A}{A}\quad\text{ and}\quad
   \w{\rho^\oplus_A}{A\oplus0}A\\
&\w{\gamma^\oplus_{A,B}}{A\oplus B}{B\oplus A}\\
&\w{\alpha^{\otimes}_{A,B,C}}{(A\otimes B)\otimes C}{A\otimes(B\otimes C)}\\
&\w{\lambda^\otimes_A}{1\otimes A}A\quad\text{ and}\quad
\w{\rho^\otimes_A}{A\otimes1}A\\
&\w{\gamma^\otimes_{A,B}}{A\otimes B}{B\otimes A}\\
&\w{\delta_{A,B,C}}{A\otimes(B\oplus C)}{(A\otimes B)\oplus(A\otimes C)} \\
&\w{\eps_A}{A\otimes 0}0
\end{align*}
be as described above in this section. Then these isomorphisms satisfy all the BD axioms.
\end{thm}

\begin{proof}
We need to prove that, with the given isomorphisms, (the diagrams in) Figures~1--18 commute. We have already observed that Figures~1--4 commute because, when category-theoretic products exist, they always provide a symmetric monoidal structure. And we have assumed a braided monoidal structure for $\otimes$, so Figures~5--8 commute.  It remains, therefore, to check the twelve requirements in Figures~9--18. (Recall that Figures~9 and 18 have two diagrams each, so these ten figures impose twelve requirements.) Fortunately, five of the twelve are trivial: Figures~14, 16, the second diagram in Figure~9 and both diagrams in Figure~18 involve morphisms to the zero object. There is just one morphism from any given object to 0, so these five diagrams automatically commute. (In fact, all the objects in these five diagrams are 0, so we could also have inferred commutativity of the diagrams from the fact that 0 is an initial object.)

So we still have seven diagrams to check for commutativity, namely Figures~10--13, 15, 17 and the first diagram in Figure~9. Recall, from our description in \S\ref{sec:axioms} of how to interpret cyclic diagrams of isomorphisms, that it suffices to check, for one pair of vertices Start and Finish, that the two paths from Start to Finish represent the same isomorphism. It then follows that the same is true for any other pair of vertices.

We will usually choose Finish to be the sum (by $\oplus$) of some objects.  Then, to check equality of two morphisms into Finish, it suffices to check equality after composing the two morphisms with the projections to the summands of Finish.  That is, if $\Pi_1$ and $\Pi_2$ are (the morphisms represented by) the two paths from Start to Finish, and if $\text{Finish}=\bigoplus_iZ_i$ with projections $p_i:\text{Finish}\to Z_i$, then to prove that $\Pi_1=\Pi_2$ it suffices to prove that $\Pi_1\circ p_i\ =\ \Pi_2\circ p_i$ for all $i$.  We call this the \emph{reduction to summands} method. We often use the same name for a path and the isomorphism represented by the path; this little abuse of notation allows us to avoid excessive pedantry.

Most of our seven proofs will begin by specifying vertices Start and Finish in the diagram, indicating the relevant paths $\Pi_1$ and $\Pi_2$, and then computing the compositions $\Pi_1\circ p_i$ and $\Pi_2\circ p_i$ for the projections $p_i$ into the summands of Finish.

Turning to the first diagram in Figure~9, let Start be the upper left corner $A\otimes(B\oplus C)$ of the diagram, and let Finish be the lower right corner $(A\otimes B)\oplus(A\otimes C)$.  Consider first the path $\Pi_1$ from Start to Finish that goes across the top of the diagram and then down the right side. When followed by the first projection, $p_{A\otimes B,A\otimes C,1}$, it gives
\begin{align*}
&\Pi_1\circ p_{A\otimes B,A\otimes C,1}\ =\
\\
&\delta_{A,B,C}\circ(\beta_{A,B}\oplus\beta_{A,C})\circ
p_{A\otimes B,A\otimes C,1}\ =\
\\
&\delta_{A,B,C}\circ p_{A\otimes B,A\otimes C,1}\circ\beta_{A,B}\ =\ \\
&(1_A\otimes p_{B,C,1})\circ\beta_{A,B},
\end{align*}
where we used the definition of $\oplus$ on morphisms and the definition of $\delta$.  Let $\Pi_2$ be the other path from Start to Finish, going down the left side of the diagram and then across the bottom. We have
\begin{align*}
&\Pi_2\circ p_{A\otimes B,A\otimes C,1}\ =\
\\
&\beta_{A,B\oplus C}\circ\delta_{A,B,C}\circ p_{A\otimes B,A\otimes C,1}\ =\
\\
&\beta_{A,B\oplus C}\circ(1_A\otimes p_{B,C,1})\ =\
\\
&(1_A\otimes p_{B,C,1})\circ\beta_{A,B},
\end{align*}
where we used the definition of $\delta$ and the naturality of $\beta$
(with respect to $1_A$ and $p_{B,C,1}$). Since the last lines in the
two computations agree, we have established that the isomorphisms
represented by the two paths have the same composite with the first
projection. The proof for the second projection is the same --- just
change 1 to 2 in the subscripts of $p$ and change $\beta_{A,B}$ to
$\beta_{A,C}$ in the preceding computation.  This completes the proof
of commutativity for Figure~9.

In Figure~10, we use the paths from the upper left to the upper right corner. Let $\Pi_1$ be the long path, and $\Pi_2$ the short one.  Composing $\Pi_1$ with the first projection gives
\begin{align*}
&\Pi_1\circ p_{A\otimes B,A\otimes C,1}\ =
\\
&(1_A\otimes\gamma^\oplus_{B,C})\circ\delta_{A,C,B}\circ
\gamma^\oplus_{A\otimes C,A\otimes B}\circ
 p_{A\otimes B,A\otimes C,1}\ =
 \\
&(1_A\otimes\gamma^\oplus_{B,C})\circ\delta_{A,C,B}\circ
p_{A\otimes C,A\otimes B,2}\ =
\\
&(1_A\otimes\gamma^\oplus_{B,C})\circ(1_A\otimes p_{C,B,2})\ =
\\
&1_A\otimes p_{B,C,1},
\end{align*}
where we used the definitions of $\gamma^\oplus$ and $\delta$, followed by the functoriality of $\otimes$ and a second use of the definition of $\gamma^\oplus$. Composing $\Pi_2$ with the first projection gives
\begin{align*}
&\Pi_2\circ p_{A\otimes B,A\otimes C,1}\ =
\\
&\delta_{A,B,C}\circ p_{A\otimes B,A\otimes C,1}\ =
\\
&1_A\otimes p_{B,C,1},
\end{align*}

Repeating the same calculation for the second projection completes the proof, by reduction to summands, of commutativity for Figure~10.

The argument for Figure~11 is more involved.
Let Start and Finish be the lower left and the upper right corners respectively. Let $\Pi_1, \Pi_2$ be the clockwise and counterclockwise paths from Start to Finish. By reduction to summands, it suffices to verify the following three equalities:
\begin{align}
\Pi_1\circ p_{AB,AC\oplus AD,1}
                   \phantom{\circ p_{AC,AD,1}}\
&=\
              \Pi_2\circ p_{AB,AC\oplus AD,1}
\tag{E11.1}\\
\Pi_1\circ p_{AB,AC\oplus AD,2}
                   \circ p_{AC,AD,1}\
&=\
              \Pi_2\circ p_{AB,AC\oplus AD,2}\circ p_{AC,AD,1}
\tag{E11.2}\\
\Pi_1\circ p_{AB,AC\oplus AD,2}
                   \circ p_{AC,AD,2}\
&=\
              \Pi_2\circ p_{AB,AC\oplus AD,2}\circ p_{AC,AD,2}
\tag{E11.3}
\end{align}

\smallskip
Figure~11.1 illustrates the proof of equality (E11.1).

\begin{figure}[H]\footnotesize
\[\xymatrix@C-2pc@R-1pc{
A(B\op (C\op D)) \ar[rr]^{\delta}
           \ar[ddrr]_{1\x p}
&&AB\op A(C\op D) \ar[rr]^{1\op \delta} \ar[dd]_p
&&AB\op (AC\op AD) \ar[ddll]^p
\\
&\mbox{} \ar @{} [ru] | {\boxed{\delta}\quad}
&&\mbox{} \ar @{} [lu] | {\quad\boxed{\op }}
\\
&&AB
\\
&\boxed{\ap}\hspace{4em}
&\boxed{\delta}
&\hspace{4em}\boxed{\ap}
\\
&A(B\op C) \ar[ruu]^{1\x p} \ar[rr]_\delta
&&AB\op AC \ar[luu]_p
\\
&\boxed{\delta}
&\hspace{6em}\boxed{\op }
\\
A((B\op C)\op D) \ar[uuuuuu]^{1\x\ap}
           \ar[ruu]^{1\x p}
           \ar[rr]_{\delta}
&& A(B\op C)\op AD \ar[luu]_p
             \ar[rr]_{\delta\x1}
&&(AB\op AC)\op AD \ar[luu]_p
             \ar[uuuuuu]_{\ap}
}\]
\caption{11.1}
\end{figure}

Let $L$ and $R$ be the left and right parts of E11.1 respectively. $L$ is the (isomorphism represented by the) extension of $\Pi_1$ with the arrow from Finish to $AB$. Similarly, $R$ is the extension of $\Pi_2$ with the same arrow. Let $M$ be the two-leg path in Figure~11.1 from Start to $A(B\op C)$ and then to $AB$.

The triangle on the left of Figure~11.1 commutes by the definitions of $\alpha^{\op}$. The two triangles in the upper part of the figure commute by the definitions of $\delta$ and $\op$ respectively. It follows that the isomorphisms represented by $L$ and $M$ coincide.

The remaining three triangles of Figure~11.1 commute by the definition of $\delta$ and $\alpha^{\op}$, and the quadrangle commutes by the definition of $\op$. It follows that the isomorphisms $R$ and $M$ coincide. Thus $L = M = R$.

In a similar way, Figures~11.2 and 11.3 illustrate the proofs of equalities (E11.2) and (E11.3).

\begin{figure}[H]\footnotesize
\[\xymatrix@C-2pc@R+0pc{
A(B\op (C\op D)) \ar[rr]^{\delta}
           \ar[dr]_{1\x p}
&\mbox{}   \ar @{} [d] | {\boxed{\delta}}
&AB\op A(C\op D) \ar[rr]^{1\op \delta} \ar[dl]^p
           \ar @{} [dr] | {\boxed{\op }}
&&AB\op (AC\op AD) \ar[dl]^p
\\
&A(C\op D) \ar[rr]^\delta
        \ar[dr]_{1\x p}
&\mbox{}\ar @{} [d] | {\boxed{\delta}}
&AC\op AD  \ar[dl]^p
\\
&{\boxed{\ap}}
&AC    \ar @{} [rr] | {\boxed{\ap}}
        \ar @{} [d]  | {\boxed{\delta}}
&&\mbox{}
\\
&A(B\op C) \ar[ru]^{1\x p}
        \ar[rr]_\delta
&\mbox{}
&AB\op AC  \ar[lu]_p
\\
A((B\op C)\op D) \ar[rr]_{\delta}
           \ar[ru]^{1\x p}
           \ar[uuuu]^{{1_A}\x\ap}
&\mbox{}   \ar @{} [u] | {\boxed{\delta}}
& A(B\op C)\op AD\ar[rr]_{\delta\x1}
            \ar[lu]_p
            \ar @{} [ru] | {\boxed{\op }}
&&(AB\op AC)\op AD\ar[uuuu]_{\ap}
            \ar[lu]^p
}\]
\caption{11.2}
\end{figure}

\begin{figure}[H]\footnotesize\label{fig:11.3}
\[\xymatrix@C-2pc@R-1pc{
A(B\op (C\op D)) \ar[rr]^{\delta}
           \ar[ddr]_{1\x p}
&&AB\op A(C\op D) \ar[rr]^{1\op \delta}
            \ar[ldd]_p
            \ar@{}[ddr] | {\boxed{\op }}
&&AB\op (AC\op AD) \ar[ddl]^p
\\
&\boxed{\delta}
\\
&A(C\op D) \ar[rr]^\delta
        \ar[ddr]_{1\x p}
&&AC\op AD \ar[ddl]^p
\\
&\boxed{\ap}\hspace{4em}
&\boxed{\delta}
&\hspace{4em}\boxed{\ap}
\\
&&AD
\\
&\mbox{}\phantom{A} \ar@{}[r]|(.65){\boxed{\delta}}
&&\mbox{}\phantom{B} \ar@{}[l]|(.65){\boxed{\op }}
\\
A((B\op C)\op D) \ar[uuuuuu]^{1\x\ap}
           \ar[rruu]^{1\x p}
           \ar[rr]_{\delta}
&& A(B\op C)\op AD \ar[uu]_p
             \ar[rr]_{\delta\x1}
&&(AB\op AC)\op AD \ar[lluu]_p
             \ar[uuuuuu]_{\ap}
}\]
\caption{11.3}
\end{figure}

For Figure~12, we consider the paths from the upper left to the lower
left corners. In contrast to the preceding computations, we work
directly with the isomorphisms in the diagram, rather than composing
them with projections.  The composite of the vertical map on the right
and the horizontal map on the bottom of the diagram is simply the
first projection $p_{AB,A0,1}:(AB)\oplus(A0)\to AB$, by the definition
of $\oplus$ of morphisms and the definition of $\rho^\oplus$ as the
first projection. Composing that with the $\delta$ across the top of
the diagram yields $1_A\otimes p_{B,0,1}$ by definition of
$\delta$. Since $p_{B,0,1}=\rho^\oplus_B$, this calculation shows that
the long path from upper left to lower left represents
$1_A\otimes\rho^\oplus_B$.  That coincides with the short path, so the
diagram commutes.

For Figure~13, we consider the paths from the upper left corner to the
lower right, and we compose them with the projection from the lower
right $A(BC)\oplus A(BD)$ to $A(BC)$.  The computation using the
projection to $A(BD)$ is exactly analogous, so we omit it.  Consider
first the path that goes across the top and then down the right side
of the diagram. Beginning with $p_{A(BC),A(BD),1}$ at the lower right
corner, we compose it with the last morphism $\delta_{A,BC,BD}$ on our
path obtaining, by definition of $\delta$, $1_A\otimes p_{BC,BD,1}$.
Composing with the next morphism back along the path,
$1_A\otimes\delta_{B,C,D}$, remembering functoriality of $\otimes$ and
the definition of $\delta$, we get $1_A\otimes(1_B\otimes p_{C,D,1})$.
Next, we must compose this with the first morphism in our path
$\alpha^\otimes_{A,B,C\oplus D}$. By naturality of $\alpha^\otimes$
(with respect to $1_A$, $1_B$, and $p_{C,D,1}$), the result is
$((1_A\otimes 1_B)\otimes p_{C,D,1})\circ\alpha^\otimes_{A,B,C}$,
which, by functoriality of $\otimes$, is
$(1_{AB}\otimes p_{C,D,1})\circ\alpha^\otimes_{A,B,C}$.  Compare this
with the other path, going down the left side and across the bottom of
the diagram.  Beginning with $p_{A(BC),A(BD),1}$ at the lower right
corner, we compose it with the last morphism
$\alpha^\otimes_{A,B,C}\oplus \alpha^\otimes_{A,B,D}$ on the path,
obtaining $p_{(AB)C,(AB)D,1}\circ\alpha_{A,B,C}$ by definition of
$\oplus$ on morphisms. Composing this with the $\delta_{AB,C,D}$ on
the left side of the diagram, and remembering the definition of
$\delta$, we get $(1_{AB}\otimes
p_{C,D,1})\circ\alpha^\otimes_{A,B,C}$, in agreement with what we
found for the first path. This completes the proof of commutativity
for Figure~13.

For Figure~15, we consider the two paths from the upper left corner
to the bottom vertex, and we compose them with the projection
$p_{A,B,1}:A\oplus B\to A$.  The computation for the other projection,
to $B$, is exactly analogous and therefore omitted.  For the longer
path, composing $p_{A,B,1}$ with
$\lambda^\otimes_A\oplus\lambda^\otimes_B$ produces
$p_{1A,1B,1}\circ\lambda^\otimes_A$ by definition of $\oplus$ on
morphisms. Composing this with the remaining morphism $\delta_{1,A,B}$
on the path produces, by definition of $\delta$, $(1_1\otimes
p_{A,B,1})\circ\lambda^\otimes_A$. By naturality of $\lambda$ (with
respect to $p_{A,B,1}$), this is the same as $\lambda^\otimes_{A\oplus
  B}\circ p_{A,B,1}$. But this is what we would get by using the
shorter path instead, so the commutativity of Figure~15 is
established.

Finally, we turn to Figure~17. Although there are 13 isomorphisms in
this diagram, many of them are essentially trivial for our
purposes. We treat the trivial ones first, namely the bottom two
morphisms in the left and right columns and the horizontal morphism
across the bottom.  These involve only additive associativity and
commutativity, so we can easily analyze the morphism represented by
the path, from left to right, consisting of these five morphisms.

In Figure~17.0, $\sigma$ is the isomorphism represented by the path in Figure~17 that goes from the third-from-bottom element on the left to the third-from-bottom element on the right. It suffices to prove that Figure~17.0 commutes.

\begin{figure}[H]\footnotesize
\[\xymatrix@C-2pc@R+0pc{
&(A\oplus B)(C\oplus D)\ar@/^/[rd]^{\delta_{A\oplus B,C,D}}
           \ar@/_/[ld]_{\gx_{A\oplus B,C\oplus D}}
\\
(C\oplus D)(A\oplus B) \ar[d]_{\delta_{C\oplus D,A,B}}
&&(A\oplus B)C\oplus(A\oplus B)D\ar[d]^{\gx_{A\oplus B,C}\oplus\gx_{A\oplus B,D}}
\\
(C\oplus D)A\oplus(C\oplus D)B\ar@{<-}[d]_{\gx_{A,C\oplus D}\oplus\gx_{B,C\oplus D}}
&&C(A\oplus B)\oplus D(A\oplus B)\ar[d]^{\delta_{C,A,B}\oplus\delta_{D,A,B}}
\\
A(C\oplus D)\oplus B(C\oplus D)\ar[d]_{\delta_{A,C,D}\oplus\delta_{B,C,D}}
&&((CA)\oplus(CB))\oplus((DA)\oplus(DB)) \ar@{<-}[d]_%
  {(\gx_{A,C}\oplus\gx_{B,C})}^{\oplus\ (\gx_{A,D}\oplus\gx_{B,D})}
\\
((AC)\oplus(AD))\oplus((BC)\oplus(BD))
\ar[rr]^{\sigma}
&&((AC)\oplus(BC))\oplus((AD)\oplus BD))
}\]
\caption{17.0}
\end{figure}

For Figure~17.0, let Start and Finish be the top vertex $(A\oplus B)(C\oplus D)$ and the bottom right vertex $((AC)\oplus(BC))\oplus((AD)\oplus(BD))$ respectively. Let $\Pi_1$ and $\Pi_2$ be the counterclockwise and the clockwise paths from Start to Finish respectively. In virtue of the reduction to summands method, it suffices to prove  the equality
\begin{equation}\label{eq:17a}
\Pi_1\circ p_{AC\oplus BC,AD\oplus BD,1} \circ p_{AC,BC,1}\
=\
\Pi_2 \circ p_{AC\oplus BC,AD\oplus BD,1} \circ p_{AC,BC,1}
\end{equation}
and similar equalities for $AD$, $BC$, or $BD$ in place of $AC$. The four proofs are similar, and so we give details only for the case of $AC$.

The path $\Pi_1$ is composed of five arrows. Let $\Pi_0$ be the initial segment of $\Pi_1$ composed of the first four arrows, so that $\Pi_1 = \Pi_0\circ \sigma$. Consider the composition of $\sigma$ with the projection to $AC$. Chasing through the definitions of $\gamma^\oplus$ and $\alpha^\oplus$, we get
\begin{equation}\label{eq:sigma}
\sigma; p_{AC\oplus BC,AD\oplus BD,1} \circ p_{AC,BC,1}\
=\ p_{AC\oplus BC,AD\oplus BD,1} \circ p_{AC,BC,1}
\end{equation}
Accordingly, \eqref{eq:17a} is equivalent to the equation
\begin{equation}\label{eq:17b}
\Pi_0\circ p_{AC\oplus BC,AD\oplus BD,1} \circ p_{AC,BC,1}\
=\
\Pi_2 \circ p_{AC\oplus BC,AD\oplus BD,1} \circ p_{AC,BC,1}
\end{equation}

To prove equation \eqref{eq:17b}, we analyze the paths in the left and right parts of \eqref{eq:17b} separately. Figure~17.1 illustrates the analysis of the left path which is the path from the top right corner to the bottom right corner of Figure~17.1. The figure splits into 3 triangles and 3 quadrangles. All 6 component diagrams commute.

\begin{figure}[H]\footnotesize\label{17.1}
\[\xymatrix@C+3pc@R+1pc{
&&(A\oplus B)(C\oplus D)
  \ar@/_/[lld]_*+{\gx}
  \ar[dddd]|(.5)*+{p_{A,B,1}\otimes p_{C,D,1}}
  \ar@/_/[lddd]_*+{p\x1}
\\
(C\oplus D)(A\oplus B)
   \ar[d]_*+{\delta}
   \ar[dr]^{1\x p}
\\
(C\oplus D)A\oplus(C\oplus D)B
  \ar@{<-}[d]_*+{\gx\oplus\gx}
  \ar[r]^p
  \ar@{}[rruu]|(.15)*+{\boxed\delta}
  \ar@{}[rruu]|(.5)*+{\boxed{\gx\text{ is natural}}}
&(C\oplus D)A
&\mbox{}
\\
A(C\oplus D)\oplus B(C\oplus D)
  \ar[d]_*+{\delta\oplus\delta}
  \ar[r]^p
  \ar@{}[ru]|(.5)*+{\boxed\oplus}
&A(C\oplus D)
  \ar[u]^*+{\gx}
  \ar[d]_*+{\delta}
  \ar[rd]^*+{1\x p}
  \ar@{}[ru]|(.55)*+{\boxed{\x\text{ is a bifunctor}}}
\\
(AC\oplus AD)\oplus(BC\oplus BD)
  \ar[r]_*+{p}
  \ar@{}[ru]|(.5)*+{\boxed\oplus}
&AC\oplus AD
  \ar@{}[ruuuu]|(.14)*+{\boxed\delta}
  \ar[r]_*{p}
&AC
}\]
\caption{17.1}
\end{figure}

The figure gives the reason for the commutativity of all six component diagrams. It follows that the perimeter of the figure commutes as well, and therefore we have
\begin{equation}\label{eq:left}
\Pi_0\circ p_{AC\oplus BC,AD\oplus BD,1} \circ p_{AC,BC,1}\
=\
p_{A,B,1}\otimes p_{C,D,1}
\end{equation}

Figure~17.2 illustrates the analysis of the right path of \eqref{eq:17b} which is the path from the top left corner to the bottom left corner of Figure~17.2.

\begin{figure}[H]\scriptsize\label{17.2}
\[\xymatrix@C+1pc@R+3pc{
(A\oplus B)(C\oplus D)
  \ar@/^3pc/[rrrd]^*+{\delta}
  \ar@{}[rrrd]|(.5){\boxed\delta}
  \ar@/_5pc/[dddd]|(.2)*+{p_{A,B,1}\x p_{C,D,1}}
  \ar[drr]_*+{1\x p}
  \ar@{}[dddd]|(.3){\boxed{\text{$\x$ is a bifunctor}}}
  \ar@{}[dddr]|(.75){\quad\boxed{\gx\text{ is natural}}}
\\
&&(A\oplus B)C
  \ar@/_7pc/[dddll]|(.5)*+{p\x1}
  \ar[d]_*+{\gx}
&(A\oplus B)C\oplus(A\oplus B)D
  \ar[d]^*+{\gx\op\gx}
  \ar[l]^(.6)*+{p}
\\
&&C(A\oplus B)
  \ar[dl]_(.7)*+{1\x p}
  \ar@{}[ru]|(.5){\boxed{\op}}
  \ar[d]^*+{\delta}
&C(A\oplus B)\oplus D(A\oplus B)
  \ar[d]^*+{\delta\op\delta}
  \ar[l]^(.6)*+{p}
\\
&CA
  \ar@{}[rruu]|(.2){\qquad\boxed\delta}
&CA\oplus CB
  \ar[l]^(.6)*+{p}
  \ar@{}[ru]|(.5){\boxed{\op}}
&(CA\oplus CB)\oplus(DA\oplus DB)
  \ar[l]^(.63)*+{p}
\\
AC
  \ar[ru]^*+{\gx}
  \ar@{}[rru]|(.5){\qquad\quad\boxed{\op}}
&&AC\oplus BC
  \ar[ll]^*+{p}
  \ar[u]|*+{\gx+\gx}
  \ar@{}[ru]|(.5){\boxed{\op}\phantom{.75em}}
&(AC\oplus BC)\oplus(AD\oplus BD)
  \ar[u]|(.5)*+{(\gx\oplus\gx) \oplus (\gx\oplus\gx)}
  \ar[l]^(.63)*+{p}
}\phantom{mmm}\]
\caption{17.2}
\end{figure}

The figure splits into 4 triangles and 4 quadrangles. All these eight component diagrams commute for the reasons indicated in the figure. It follows that the perimeter of the figure commutes as well, and therefore we have
\begin{equation}\label{eq:right}
\Pi_2 \circ p_{AC\oplus BC,AD\oplus BD,1} \circ p_{AC,BC,1}\
=\
p_{A,B,1}\otimes p_{C,D,1}
\end{equation}

Equation \eqref{eq:17b} follows from \eqref{eq:left} and \eqref{eq:right}. This completes the proof of commutativity for Figure~17 and thus completes the proof of the theorem.
\end{proof}

\end{document}